\documentclass[12pt,notitlepage,twoside,a4paper]{amsart}
 \usepackage{amsfonts}

\usepackage{amsmath,amssymb,enumerate}

\usepackage[usenames, dvipsnames]{color}

\usepackage{xcolor}
\usepackage{lipsum}
\usepackage{amssymb}
\usepackage{amsmath,amsthm}
\usepackage{latexsym}
\usepackage{amscd}
\usepackage{psfrag}
\usepackage{graphicx}
\usepackage[latin1]{inputenc}
\usepackage[all]{xy}
\usepackage[mathcal]{eucal}

\definecolor{NoteColor}{rgb}{1,0,0}

\renewcommand{\textsc}{\textcolor{red}}

\makeatletter
\newcommand\Max{\@tempcnta=\mathcode`\m\relax
\mathcode`\m=\mathcode`\M\max\mathcode`\m=\@tempcnta\relax}

\newcommand\Min{\@tempcnta=\mathcode`\m\relax
\mathcode`\m=\mathcode`\M\min\mathcode`\m=\@tempcnta\relax}
\makeatother

\newtheorem*{theorem 1}{\rm\bf Proposition 1}
\newtheorem*{theorem 2}{\rm\bf Proposition 2}

\newcommand{\M}{\mathfrak{M}}

\theoremstyle{definition}

\theoremstyle{remark}

\def\interieur#1{\mathord{\mathop{\kern 0pt #1}\limits^\circ}}

\title[Value distribution theory]{Value distribution theory and Teichm\"uller's paper \emph{Einfache Beispiele zur Wertverteilungslehre}}

\author{Athanase Papadopoulos}

\address{Athanase Papadopoulos, Institut de Recherche Mathématique Avancée, Universit\'{e} de Strasbourg et CNRS, 7 rue Ren\'{e} Descartes, 67084 Strasbourg Cedex, France 
\\
 {\tt papadop@math.unistra.fr} }

\thanks{}

\date{\today}

\makeindex

\begin{document}

  \maketitle
  \begin{abstract}
   
This survey will appear in Vol. VII of the \emph{Hendbook of Teichm\"uller theory}. It is a commentary on Teichm\"uller's paper  ``Einfache Beispiele zur Wertverteilungslehre", published in 1944,  whose English  translation appears in that volume.  Together with Teichm\"uller's paper, we survey the development of value distribution theory, in the period starting from Gauss's work on the Fundamental Theorem of Algebra and ending with the work of Teichm\"uller. We mention the foundational work of  several mathematicians, including Picard, Laguerre, Poincaré, Hadamard, Borel, Montel, Valiron, and others, and we give a quick overview of the various notions introduced by Nevanlinna and some of his results on that theory.

\end{abstract}

          AMS classification: 30D35, 30D30
          
          Keywords: value distribution theory, Picard theorems, Nevanlinna theory, Nevanlinna theorems, order of a meromorphic function, proximity function, Nevanlinna characteristic function, exceptional values of an entire function, deficiency function (Nevanlinna), exceptional value of a meromorphic function.

  \tableofcontents
\section{Introduction}

Value distribution theory\index{value distribution theory} is the study of the distribution of the values taken by a meromorphic function. It originates in works of Weierstrass, Laguerre, Hadamard and Poincaré, and especially  in Picard's famous two theorems and the attempts by several mathematicians to give more precise statements, simpler proofs and generalizations of them. The theory experienced a spectacular progress in the work done by Nevanlinna\index{Nevanlinna, Rolf (1895--1980)} in the years 1922--1925 and for this reason it is also called Nevanlinna's theory. The references to Nevanlinna's work of that period are the papers \cite{N1922, N1924, N1924a, N1924b, N1924c, N4}.
Riemann surfaces were introduced in the theory, through the notion of Riemann surface associated to a meromorphic function. This  brought into value distribution theory topological ideas, especially those related to the techniques of ramified coverings of the sphere. Methods from differential geometry and quasiconformal mappings were also introduced. After Nevanlinna, Ahlfors (who was Nevanlinna's first student) and Teichm\"uller were among the main promoters of that new turn.\index{value distribution theory}

In this chapter, we survey some of the important ideas and developments of value distribution theory\index{value distribution theory} from the early period until the work of Teichm\"uller.  

Teichm\"uller probably became familiar with value distribution theory by reading Nevanlinna's book \emph{Le th\'eor\`eme de Picard--Borel et la th\'eorie des fonctions m\'eromorphes} \cite{N1} published in 1929, translated into German in 1936, and which he quotes in his paper \cite{T33} published in 1944. He also attended  Nevanlinna's lectures in G\"ottingen in 1936 as well as lectures by Egon Ullrich\index{Ullrich, Egon (1902--1957)} on this topic.  In his paper \emph{Extremale quasikonforme Abbildungen und quadratische Differentiale} \cite{T20}, published in 1939, he\index{Teichm\"uller, Oswald (1913--1943)} mentions the importance of  quasiconformal mappings in value distribution theory which motivated several of his works, including the papers \cite{T13} (1938) and \cite{T33} (1944) translated in the present volume.

Between the works of Gauss on the Fundamental Theorem of Algebra and those of Teichm\"uller, we shall review some foundational works of Picard,  Laguerre, Hadamard, Borel, Montel, Valiron, and others, and we shall  highlight the main notions introduced by Nevanlinna in this field.

  A very good survey on the early Nevanlinna theory is Lehto's paper \emph{On the birth of the Nevanlinna theory}  \cite{Lehto-birth}. 
  
  I would like to thank Vincent Alberge for his remarks on a first version of this paper.
  
   \section{Value distribution theory before Nevanlinna}\label{s:val}
  
Let $f$ be a meromorphic function defined on the complex plane and let  $a$ be a complex number. One would like to have precise information about the distribution of the solutions of the equation 
$f(z)=a$, counted with or without multiplicity. The kind of information includes an estimate of  
the number of solutions in a disc $\{\vert z\vert \leq r\}$ for any positive real number $r$, estimates on the growth and the asymptotics of the number of such solutions in terms of $r$, the comparison of the various estimates when the constant $a$ varies, etc. Value distribution theory \index{Nevanlinna theory} was born when mathematicians started asking such questions. 

           A notable ancestor of this theory is the Fundamental Theorem of Algebra,\index{Fundamental Theorem of Algebra}\index{Theorem!Fundamental Theorem of Algebra} stating that if $f$ is a polynomial, then, for each complex number $a$, the equation $f(z)=a$  has at least one solution. Several attempts of proofs of this theorem were given in the eighteenth century by d'Alembert, Euler and others, but they are all considered as incomplete, as Gauss reports in his doctoral dissertation in 1799 \cite{Gauss1799}, in which he gives a proof of that theorem. It turned out that Gauss's proof contains a gap of a  topological nature, which was filled by Ostrowski  in 1920 \cite{Ostrovski1920}. In the meanwhile, several other proofs of the Fundamental Theorem of Algebra were discovered, including seven different proofs by Gauss himself. For the history of this problem, we refer the reader to the survey article by Remmert \cite{Remmert1991}.

     Weierstrass began a systematic investigation of the theory of zeros of entire functions. He proved in 1876 that for any sequence  $(a_n)$ of complex numbers whose moduli $\vert a_n\vert$ are increasing and  tend to infinity,  one can find a holomorphic function defined on the plane and having the set $\{a_n\}$ as zero set \cite{W34}. He constructed such a function as an infinite product. This result is a wide counterpart of the fact that given a finite sequence of complex numbers, one can write (using a product formula) a polynomial having these numbers as zeros.  Weierstrass's result is  one of the first results on the natural  attempt to generalize to entire functions results known for solutions of  polynomial equations.

           An important step towards value distribution theory originates in the desire to generalize notions like the growth of 
           a polynomial and the properties of that growth to more general analytic functions.
           
           To be more precise, let us consider again a polynomial $f$ of degree $n\geq 1$ and let $a$ be and arbitrary complex number. Then, the number of points  $z$ in the complex plane satisfying $f(z)=a$, counted with multiplicity, is exactly $n$. Furthermore, the growth of the value $f(z)$ when  $z\to\infty$ is also equal to $n$; more precisely,       \[\lim _{x\to \infty} \frac{f(z)}{z^n}=A,\]
           where $A$ is the coefficient of the monomial of highest degree of $f$. Thus, the (polymonial) growth of $f$ can be read on its coefficients.
           
            Value distribution theory is concerned with the search for analogous results for functions that are more general than polynomials.

                                        Laguerre\index{Laguerre, Edmond (1834--1886)} was one of the other early promoters of such a theory. In
   the introduction of his paper  \cite{Laguerre-Equations} (1882), he writes: ``The theorems of Rolle and Descartes\footnote{These are results of Rolle and Descartes concerning  solutions of algebraic equations which Laguerre discussed in previous papers, cf. \cite{L1}   and \cite{L2}.} apply to transcendental functions.  But this is not the case for the 
           consequences, which are so simple and so numerous, which we deduce from these two propositions, relatively to the equations whose left-hand side is an entire polynomial; they subsist only exceptionally. [...] Indeed the transcendental functions $\cos x$ and $\sin x$ share the properties of entire polynomials; but this is no more the case for the holomorphic function $G(x)$, the inverse of Legendre's function $\Gamma(x)$ which was introduced in analysis by Weierstrass. [...] Thus, it seems of some interest to study what are the elementary properties of algebraic equations that apply to transcendental equations."

     It is generally considered that modern value distribution theory\index{value distribution theory} was born with the two theorems of Picard\index{Picard theorem}\index{Theorem!Picard}\index{Picard, Emile (1856--1941)} which we recall now.
           
The first Picard theorem (also called \emph{Picard's small theorem})\index{Picard theorem (small)}\index{Theorem!Picard (small)} is a generalization of Gauss's theorem. It states that for any non-constant entire function $f$,  the equation $f(z)=a$ has a solution for every value of $a$ except possibly for one value. (An example of an entire function with an exceptional value is the  exponential map $f(z)=e^z$, with $a=0$.) 
Picard obtained that theorem in 1879 \cite{Picard1879}. In the same year, he  proved his second theorem, namely, that an arbitrary meromorphic function, in the neighborhood of an essential singularity, takes any complex value infinitely often, with at most one exception  \cite{Picard1879a}. This is usually called \emph{Picard's big theorem}\index{Picard theorem (big)}\index{Theorem!Picard (big)}. The small theorem follows from the big one since an entire function is either a polynomial or has an essential singularity at infinity.

Picard's big theorem is an important improvement of the result of Weierstrass\index{Weierstrass, Karl (1815--1895)} that we mentioned. One should also mention here the so-called Weierstass--Casorati\index{Weierstass--Casorati theorem}\index{Theorem!Weierstass--Casorati} theorem stating that if $f$ is a holomorphic function defined on a punctured disc of radius $r>0$ with an essential singularity at the puncture, then for every $0<r'<r$ the image of the sub-punctured disc of radius $r'$ is dense in the complex plane.

 It was considered as a surprise that Picard, in the proof of his theorems, did not use Weierstrass's theory. His proof was geometric, and one of the ingredients was the so-called \emph{elliptic modular function}, a function defined on the upper half-plane, which is invariant under the action of the modular group $\mathrm{PSL}(2,\mathbb{Z})$ and which he had already used in his previous works. One proof of Picard's small theorem goes as follows: If $f:\mathbb{C}\to \mathbb{C} $ is an entire non-constant function that omits two values, then, composing $f$ with the inverse of the modular function provides a non-constant map from the complex plane onto the upper-half plane. But the latter is biholomorphically equivalent  to the unit disc. This contradicts Liouville's theorem, saying that a bounded entire function is constant. This also reminds us of the fact that Liouville's theorem\index{Liouville theorem}\index{Theorem!Liouville} is another result which is at the foundations of function theory.  The fundamental theorem of algebra follows from it. 
 
Picard's proof of his small theorem was considered difficult, because it uses geometry rather than analysis. Simpler proofs\footnote{The notion of ``simple proof" is relative. Ahlfors, in his book \emph{Conformal invariants} writes: ``Our proofs of the Picard theorems is elementary not only because it avoids the modular function, but also because it does not use the monodromy theorem" \cite[p. 21]{Ahlfors-CI}.} were sought for, and value distribution\index{value distribution theory} theory is rooted in these attempts and in attempts to generalize the theorem. We mention some of them.

 Borel reduced Picard's theorem to a statement on linear relations between exponentials of entire functions. In 1896 \cite{Borel-Picard} he gave an ``elementary" proof of the small theorem, which he expanded in his later articles and books, \cite{Borel1897, Borel1903, Borel1921}.  Borel's method is based on results of Hadamard and Picard on the minimum modulus of an entire function and on the growth of its derivative.  Schottky in 1904 \cite{Schottky-Picard} used arguments similar to Borel's to prove other results of Picard. An elementary proof of Picard's small theorem that is often quoted is the one of Montel \cite{Montel1912}, which uses Schottky's proof and normal families.
   There are also proofs by Bloch \cite{Bloch-Picard}, Carath\'eodory \cite{Caratheodory-Picard}, Landau \cite{Landau1904}, Lindel\"of \cite{Lindelof-Picard}, Milloux \cite{Milloux-Picard}, Valiron \cite{Valiron-Picard} and others.  Borel's extension of Picard's theorem gives information on the values and on the distributions of exceptional values of entire or meromorphic functions, and it became known as the Picard--Borel Theorem. We shall say a few words on Borel's extension below.  The number of re-formulations, sharpening, extensions and proofs of Picard's theorems that were published is an indication of the deepness of the theorem and its relation with many parts of mathematics.
   
   Nevanlinna\index{Nevanlinna, Rolf (1895--1980)} in his book \cite{N1} reviews some proofs and generalizations of Picard's theorem.   A proof due to  Ahlfors \cite{Ahlfors-Zur} has a topological flavor. In fact, Ahlfors obtained a geometric form of Picard's theorem and of a generalization by Nevanlinna which has little to do with analytic functions. We shall quote it in \S \ref{s:Ahlfors} below.  

Let us mention some of the new notions regarding entire functions  that were introduced during these investigations.
                      
                         Using Weierstrass's work, 
            Laguerre\index{Laguerre, Edmond (1834--1886)} defined a notion of \emph{genus} of an entire function which has properties which makes it analogous to the degree of a polynomial. In his paper \cite{Laguerre-Genre} \emph{Sur la d\'etermination du genre d'une fonction transcendante enti\`ere} (1882) he gave the following characterization of the genus: 
         if $n$ is an integer such that $\frac{f'(z)}{f(z)z^n}\to 0$ as $z$ tends to infinity, then the function $f(z)$ is of genus $n$.

                      Poincar\'e\index{Poincaré, Henri (1854--1912)} in his paper \cite{P1883} \emph{Sur les fonctions enti\`eres} (1883) proved that for any entire function $f$, if \[
                      M(r)=\max_{\vert z\vert \leq r}\vert f(z)\vert
                     \]                       denotes the maximum modulus of $f$ and           
 $k$ its genus, then
\[\log M(r)=o(r^{k+1})\] as $r\to\infty.$
 Hadamard\index{Hadamard, Jacques (1865--1963)},  in his paper   \emph{\'Etude sur les propri\'et\'es des fonctions enti\`eres et en particuler d'une fonction consid\'er\'ee par Riemann}
 \cite{H1893} published in 1893, completed Poincaré's study of the notion of genus of an entire function

These results of Laguerre, Poincaré and Hadamard were obtained soon after Picard published his theorems, but the latter remained rather isolated for several years, and the activity around it started more than 15 years after its discovery. Bloch writes in \cite{Bloch1926}, in a historical survey of the theory of meromorphic functions, that there are very few examples in the history of science of a theorem of such importance which remained isolated for a long period of time.

            Borel, in 1897, defined the notion of \emph{order} \index{order (entire function)} of an entire function \cite{Borel1897}, which is a measure of how fast the maximum modulus of the function grows. More precisely, 
            he defined the order of $f$ to be the quantity
            \[\rho=\inf\{ \lambda \mid \log M(r)=O(r^{\lambda})\}.\]

             With such a definition, the order of a polynomial is zero,   the exponential function $e^z$ has order one and the function $e^{e^z}$ has infinite order. Borel's definition of order was used later by Nevanlinna who generalized it to the setting of meromorphic functions that are not necessarily entire (see below).
  
   Hadamard proved the following theorem in his 1893 paper that we quoted  \cite{H1893}: If the order $\rho$ of an entire function $f$ is finite, then for every $\epsilon >0$, the series  $\sum\vert a_n\vert^{-(\rho + \epsilon)}$ converges.
       
      Borel reformulated the Picard theorem using the notion of order and he obtained the following sharper result, called the Picard--Borel theorem:\index{Theorem!Picard--Borel}\index{Picard--Borel theorem}
          If $f$ is an entire function of order $\rho$, then we have, for any finite complex number $a$,   
          \[\lim\sup_{r\to\infty} \frac{\log n(r,a)}{\log r}\leq\rho\]
          with equality  holding except possibly for one value $a$. Here, $n(r,a,f)$ is the number of roots, with multiplicity, of the equation $f(z)=a$, in the disc $\vert z\vert \leq r$.   
           
            This is a sharpening of part of Picard's result, since the latter says that $n(r,a)$ can vanish identically for only one value $a$. 
                 
                  Borel, in his generalization of Picard's theorem,  made at the same time the connexion with the existing results on the value distribution of entire functions. 
      This generalization implies that except possibly for one value of $a$, the number of solutions of $f(z)-a=0$ which are in the disc $\vert z\vert <r$ is of the order of the logarithm of the maximum value of $\vert z\vert$ on the circle  $\vert z\vert =r$ .

          In the first years of the twentieth century, value distribution theory was already the dominant topic in the research concerning entire functions.  In 1913, Valiron published a long memoir titled \emph{Sur les fonctions entières d'ordre nul et d'ordre fini et en particulier
les fonctions à correspondance régulière} \cite{Valiron1913}. He starts his memoir by saying that research on entire functions took a new direction after the works of Hadamard  \cite{H1893} and Borel \cite{Borel1897}. He then  gives a list of mathematicians working on the questions, among them those we mentioned, but also others such as Lindel\"of (who become soon Nevanlinna's teacher), Boutroux, Blumenthal, Denjoy, Wiman, Littlewood, Sire, Mattson and Maillet.

          Julia, in 1919, gave another improvement of Picard's theorem, by showing that there exists in the complex plane a sector of arbitrarily small angle in which the equation  $f(z)-a=0$ has an infinite number of solutions, with the possible exception of one value of $a$.
          
              Valiron, at the Strasbourg ICM of 1920, addressed the question of whether one can improve Julia's result in the same way as Borel's result improves the Picard theorem.

  Let me end this section by paraphrasing a text by Painlevé, addressing Picard\index{Picard, Emile (1856--1941)} and remembering the times when he was his student:
  \begin{quote}\small
  This was the time where yourself, by a surprising synthetic effort, had just pulled off from the Unknown those two theorems on analytic functions to which your name will be attached for ever. These are revealing theorems: like two capes of an unknown continent, discovered by some daring sailors, made us anticipate a mysterious world, a world so wide and rich that fifty years of explorations did not suffice to exhaust its secrets.\footnote{C'était le temps où vous-même, par une surprenant effort synthétique, vous veniez d'arracher à l'inconnu ces deux théorèmes sur les fonctions analytiques auxquels votre nom restera à jamais attaché, théorèmes révélateurs : tels deux caps d'un continent inconnu, découverts par quelques hardis navigateurs, font pressentir un monde mystérieux, monde si vaste et si riche que cinquante années d'exploration n'en ont pas encore épuisé les secrets. (quoted by Louis de Broglie in \cite {Broglie}, in his inauguration lecture at the Académie française; de Broglie was elected at the place which was left vacant by Picard).}
  \end{quote}

 \section{Nevanlinna's two theorems, and the inverse problem} \label{s:N}
 Nevanlinna's results and methods are wide generalizations of those that Borel used to prove Picard's theorem. We shall mostly use the notation Nevanlinna sets up in his book \emph{Le th\'eor\`eme de Picard--Borel et la th\'eorie des fonctions m\'eromorphes} (1929) \cite{N1}, since this is the book that was used by Teichm\"uller. 
 
We start with a meromorphic non-rational function $f$ defined on a domain $\vert z \vert <R$ where $R$ may be any real number satisfying $0<R\leq \infty$. We recall that the aim of the theory is to give information on the distribution of the solutions of the equation 
\begin{equation}\label{eq:a}
f(z)=a
\end{equation}
where $a$ is a point in the extended complex plane $\mathbb{C}\cup\{\infty\}$,
and to compare these distributions for different values of $a$.
 
For that purpose, Nevanlinna used several functions that describe the asymptotic behavior of $f$, and we now recall some of them.

 For $r\geq 0$, the function, 
 $m(r,a)$ measures the average closeness of the value taken by the function $f$ to the complex number $a$ in the disc $\vert z\vert <r$. To make a precise definition, one starts by defining $n(r,a)$ as the number of solutions of Equation (\ref{eq:a}) in the closed disc $\vert z\vert\leq r$, counted with multiplicities.  (In the particular case where $a=\infty$, we are counting the number of poles of $f(z)$ in the given disc.)

A first \emph{counting function}\index{counting function} $N(r,a)$ is then defined as
\[N(r,a)=\int_0^r \frac{n(t,a)-n(0,a)}{t}dt +n(0,a)\log r.
\]
This function is an increasing and convex function of $\log r$. It is determined by the distribution of the moduli of the solutions of Equation (\ref{eq:a}). It measures the density of these solutions near the point $\infty$.

 We set
\[N(r,f)=N(r,\infty).\]
 
 The \emph{Nevanlinna counting function},\footnote{\label{f:1}We are using Nevanlinna's monograph \cite{N1} ; the terminology changed later.}\index{counting function (Nevanlinna)}\index{Nevanlinna counting function} also called ``smoothed counting function", is given by
\[m(r,a)= \frac{1}{2\pi}\int_0^{2\pi}\log^+\big\vert \frac{1}{f(re^{i\phi})-a}\big\vert d\phi,\]
where for each $x>0$, the \emph{positive logarithm} $\log^+ x$  is defied by 
\[\log^+ x= \max \{\log x, 0\}.\]
       (Taking positive logarithm eliminates cancellations.)

The value $m(r,a)$ is the mean value of the positive logarithm $\log^+\big\vert \frac{1}{f-a} \big\vert$ on the circle $\vert z\vert=a$ and it gives an indication of the strength of the mean convergence of the function $f(z)$ to the value $a$ as the radius $r$ tends to infinity.

In general, the function $m(r,a)$, unlike the counting function $N(r,a)$, is neither increasing nor convex in $\log r$.

At first sight, the importance of all these functions is far from being obvious. The definitions were formulated gradually in Nevanlinna's 1922--1925 papers. These papers consist of short notes \cite{N1922, N1924a, N1924c} which essentially contain announcements of results, and long papers where the results are proved in detail \cite{N1924, N1924b, N4}. 
The titles of some of these papers are informative: \emph{Sur les relations qui existent entre l'ordre de croissance d'une fonction monogène et la densité de ses zéro} (On the relations that exist between the order of growth of a monogenic\footnote{This is one of the names used for a holomorphic function (``gen" comes from a Greek work which means offspring, and it stands for the derivative; hence, the word ``monogenic" expresses the fact that a holomorphic function has a one derivative, that is, not depending on the direction).} function and the density of its zeros) \cite{N1922} (1922), \emph{Untersuchungen \"uber den Picard'schen Satz} (Researches on Picard's theorem) \cite{N1924} (1924), etc. In the paper \cite{N1922}, concerning holomorphic functions defined in an angular sector, the function $\log\vert f\vert $, which is harmonic and becomes infinite at the zeros and poles of the function $f$,  already plays an important role. From that point on, potential theory (Green's formulae, the Poisson integral, etc.) became an essential part of the theory.

Nevanlinna's monograph \cite{N1} contains a comprehensive treatment of the counting functions and related functions. Nevanlinna starts by proving several properties of the sum
\[m(r,a)+N(r,a)\] 
which turns our to be more important than the simple counting function $N(r,a)$.

He shows in particular that for a fixed value of $a$, this sum is   increasing and a convex function of $\log r$, and he studies the dependence of this function on the complex number $a$.

He defines the \emph{proximity function}\footnote{See Footnote \ref{f:1}.}\index{proximity function (Nevanlinna)}\index{Nevanlinna proximity function} as
\[m(r,f)=m(r,\infty)=\frac{1}{2\pi}\int_0^{2\pi}\log^+\vert f(re^{i\phi})\vert d\phi.\]
This function measures how large $f$ is, on the average, on the disc of center $0$ and radius $r$.

The \emph{Nevanlinna characteristic function}\index{characteristic  (Nevanlinna)}\index{Nevanlinna characteristic function} of $f$ is defined as the sum of the proximity function and the Nevanlinna counting function:
\[T(r)=T(r,f)=m(r,f) +N(r,f).\]

\emph{Nevanlinna's First Fundamental Theorem}\index{Nevanlinna First Fundamental Theorem}\index{Theorem!First Fundamental Theorem (Nevanlinna)} says that for every complex number $a$, we have
\[T(r,f)=m(r,a)+N(r,a)+h(r,a)\]
where $h(r,a)$ is a function of $r$ which stays bounded as $r\to\infty$; cf. \cite[p. 12]{N1}. 

Thus, the theorem says that if the counting function $N(r,a)$ is small (respectively large), then this is compensated by the largeness (respectively smallness) of the smooth counting function $m(r,a)$ which takes into account the strength of the convergence of $f$ to the value $a$.

The theorem was first stated in the Comptes Rendus note \cite{N1924c}. Nevanlinna later called it the ``First Fundamental Theorem."

  The first corollary that Nevanlinna deduced from his First Fundamental Theorem is that if for some value of $a$ the sum $m(r,a)+N(r,a)$ stays bounded as $r\to\infty$, then the meromorphic function $f$ is constant.

Another corollary is formulated in terms of \emph{exceptional values} of a meromorphic function $f$, that is, complex numbers $a$ for which the set of solutions of the equation $f(z)=a$ is relatively rare, or for which the growth of the function $N(r,a)$ is exceptionally slow. The corollary says that if one takes into account not only the distribution of the solutions of the equation $f(z)=a$ but also the \emph{intensity}, measured by the expression $m(r,a)$, then there are no exceptional values.

In other words, if $a$ is an exceptional value from the point of view of the distribution of the solutions of the equation $f(z)=a$, then the function $f$ converges rapidly to this value for $r\to\infty$. In contrast, if the convergence to a particular value $z$ is relatively slow, then this fact is compensated by the existence of a highly dense set of points where the value is effectively taken by the function. In either case, the sum $m+N$ is independent of the value of $a$.

Nevanlinna calls this a \emph{symmetry property} in the asymptotic shape of meromorphic functions. He declares that this property constitutes the foundations of the theory \cite[p. 13]{N1}.

 Using his counting function, Nevanlinna  defined the notion of \emph{order} of a meromorphic function by the formula
 \[\overline{\lim}_{r\to\infty}\frac{T(r)}{\log r}.\]
                This notion generalizes Borel's definition of the order of an entire function. It is used in  \emph{Nevanlinna's Second Fundamental Theorem}\index{Nevanlinna Second Fundamental Theorem}\index{Theorem!Second Fundamental Theorem (Nevanlinna)} which we now state \cite[p. 69]{N1}: 
        
        \emph{Let $f$ be a meromorphic function defined on the complex plane, let $T(r)$ be its characteristic function and let $a_1,\ldots,a_q$ ($q\geq 3$) be a finite set of distinct elements in $\mathbb{C}\cup \infty$.
        Then, 
        \[(q-2)T(r)< \sum_1^q N(r,a_q)-N_1(r)+S(r)\]
        where \[N_1(r)=N(r,\frac{1}{f'})+ (2N(r,f)-N(r,f')\]
        and where the function $S$ satisfies the following two properties:
        \begin{enumerate}
        \item For any $\lambda>0$, we have 
        \[\int_{r_0}^r \frac{S(t)}{t^{\lambda+1}}dt + \mathrm{O}\left(\int_{r_0}^r\frac{\log T(t)}{t^{\lambda+}}dt \right);\]
        \item \[S(r)<\mathrm{O}(\log T(r)+\log r),\] except possibly, in the case of infinie order, a sequence of intervals whose total length is finite.
        \end{enumerate}}

The \emph{deficiency function}\index{deficiency function (Nevanlinna)}\index{Nevanlinna deficiency function} $\delta (a)$ is defined as

 \[\delta(a)=   \liminf_{r\to R} \frac{m(r,a)}{T(r)}
 = 1- \limsup_{r\to R} \frac{N(r,a)}{T(r)}.\]
  
  We have 
  $0\leq \delta(a)\leq 1$.
  
  The value $a$ is said to be \emph{deficient} if $\delta(a)>0$. 
  
  The first results using this function obtained in the setting of Nevanlinna's theory are sufficient conditions for the equation $\delta(a)=0$ to hold.

 After proving his two theorems together with other results that appear in their developments, Nevanlinna addressed the problem which became known as Nevanlinna's inverse problem:\index{inverse problem (Nevanlinna)}   
   
  \emph{Given a (finite of infinite) sequence of numbers $\delta_n$ satisfying $0<\delta_n\leq 1$ and 
  $\sum_{n}\delta_n=2$, to find a meromorphic function admitting $n$ exceptional values $z_n$ such that
  \[\delta(z_n)=\delta_n \hbox { for all } n=1,2, \ldots. \]
  }
  
  This was formulated later in a slightly different manner, as the problem of constructing a meromorphic function with prescribed bahavior, with respect to ramification and deficiency, and  satisfying the two fundamental theorems of Nevanlinna. The first version of this problem was proposed by Nevanlinna in his book \emph{Le th\'eor\`eme de Picard--Borel et la th\'eorie des fonctions m\'eromorphes} \cite[p. 90]{N1}. Nevanlinna himself contributed to the solution  a few years later, in \cite{N2}, after he introduced a class of Riemann surfaces with finitely many logarithmic branch points. This is also one of the problems on which Teichm\"uller worked. In his paper \cite{T13}, the latter found a general principle which had several applications in value distribution theory; see \S \ref{s:Ahlfors} below.  The solution of Nevanlinna's inverse problem was completed  by Drasin in 1977 \cite{Drasin1977}.    

A more general problem, which we call the generalized inverse problem of Nevanlinna,\index{inverse problem (Nevanlinna)!generalized} is the problem of finding or characterizing functions with specific conditions on the counting function, characteristic function, etc. 
  
  Nevanlinna introduced most of the above functions in his paper \emph{Zur Theorie der meromorphen Funktionen} published in 1925 \cite{N4}. 
       Hermann Weyl,\index{Weyl, Hermann (1885--1955)} in his 1943 paper \cite{Weyl1943}, writes about that paper \cite{N4}: ``The appearance of this paper has been one of the few great mathematical events in our century." 
       
     It is (historically) important to note that at the time Nevanlinna published his work, several mathematicians were investigating similar ideas. For instance, Valiron,\index{Valiron, Georges (1884--1955)} in his paper \cite{V1925} published in 1925,  introduced a related notion of deficiency function of $a$, defined by the formula:
 \[\Delta(a)= \limsup_{r\to\infty} \frac{m(r,a)}{T(r)}= 1- \liminf_{r\to\infty} \frac{N(r,a)}{T(r)}.\]
 
 The function $\Delta(a)$ is called the \emph{Valiron deficiency}\index{Valiron deficiency function}\index{deficiency function (Valiron)} or \emph{upper deficiency}\index{upper deficiency} of $f$, cf. \cite{Ullrich1939a}. We have: 
 \[0\leq \delta(a)\leq \Delta(a)\leq 1.\]
 
 A value $a$ for which $\Delta(a)>0$ is called an \emph{exceptional value} of $f$.
  
\section{Riemann surfaces}
            
Nevanlinna, in his monograph \cite{N1}, writes (p. 289) that it is important to find the implications of his two theorems in concrete examples. As particularly interesting examples, he mentions the automorphic (Fuchsian) functions, whose associated Riemann surfaces are regularly ramified over a finite number of branch values. He then says that more generally, one may study surfaces that are branched, in a non-necessarily regular way, over finitely many values. This is how he introduced his surfaces of type $F_q$ on which we report in the paper \cite{ABP1} in the present volume.

On p. 302 
of \cite{N1}, Nevanlinna comments on a partial solution of the inverse problem\index{inverse problem (Nevanlinna)} of value distribution for surfaces with logarithmic ends, which we reformulate as follows:

\emph{Given $q$ points $w_1,\ldots,w_q$ and $q$ numbers $\delta_1,\ldots,\delta_q$ in the interval $(0,1)$ whose sum is equal to 2, to construct a meromorphic function having the $\delta_i$ as deficiency at the $w_i$ ($i=1,\ldots, q$)}.

 He mentions that the problem has been extended in an interesting way by Ullrich in his article \cite{Ullrich} and that surfaces with finitely many periodic (instead of logarithmic) ends play an important role in this work.

In Nevanlinna's 1929--1930 papers \cite{N1930} and \cite{N1930b},  Riemann surfaces associated to meromorphic functions became an important tool in value distribution theory.  In his paper \cite{N2}, he obtained a first partial solution of his inverse problem by studying a class of surfaces with finitely many logarithmic branch points. This class of surfaces played an essential role in his later work, and in the work of Teichm\"uller on value distribution theory and on the type problem; see the papers \cite{T9, T13} translated in the present volume.

\section{On Ahlfors' and Teichm\"uller's approaches} \label{s:Ahlfors}

Ahlfors introduced new geometric methods in value distribution theory, especially by developing his techniques of covering surfaces.  In some sense, this was a return to Picard's original geometric methods.  His 1935 paper \emph{Zur Theorie der \"Uberlagerungsfl\"achen} (On the theory of covering surfaces)\cite{Ahlfors-Zur} contains new proofs and generalizations of Nevanlinna's theorems, reformulated in more topological and geometric terms. One of his results says that if $f$ is a nonconstant entire function from the $z$-plane to the $w$-plane, then for any two disjoint discs in the $w$-plane, the interior of at least one of them is the image by $f$ of some domain in the $z$-plane. This is a generalization of Picard's small theorem with a topological flavor. Ahlfors also proved that (with the same notation) for any three disjoint discs in the $w$-plane, the interior of at least one of them is the one-to-one conformal image by $f$ of a domain in the $z$-plane. In the same paper, he  proved that 
if $f$ is a nonconstant meromorphic function from the $z$-plane to the $w$-Riemann sphere,  then the above conclusion holds for 5 disjoint discs instead of three. The result became known as Ahlfors' ``five-island" theorem.\index{Ahlfors' five island theorem}\index{Theorem!Ahlfors' five island} It generalizes a theorem of Nevanlinna on the number of totally ramified values of a meromorphic function of the plane. The result was conjectured by Bloch.

Perhaps more importantly, Ahlfors showed in the same paper that these results hold for quasiconformal mappings. This showed at the same time that the results had little to do with the analyticity of the maps.  He writes in the introduction to his paper \cite{Ahlfors-Zur} (translation in \cite{Lehto-Ahlfors} p. 28): 
\begin{quote}\small
This work had its origin in my endeavor to get by geometric means the most significant results of meromorphic function theory. In these attempts, it became evident that only an easily limited portion of R. Nevanlinna's Main Theorems, and thereby nearly all classical results, were dependent on the analyticity of the mapping. In contrast, their entire structure is determined by the metric and topological properties of the Riemann surface, which is the image of the complex plane. The image surface is then thought of being spread out above the Riemann sphere, i.e. to be the covering surface of a closed surface.
\end{quote}

 The topological techniques of branched coverings, cutting and pasting of pieces of the complex plane,  the Euler characteristic, the Riemann--Hurwitz formula, and other topological tools, became with Ahlfors part of the theory of functions of a complex variable. 
 Carath\'eodory, in presenting Ahlfors' work at the Fields medal ceremony in 1936, declared that the latter opened up a new chapter in analysis which could be called ``metric topology." In some sense, Ahlfors' new approach constituted a return to the sources, that is, to the geometrical and topological methods introduced by Riemann, for whom meromorphic functions were  nothing else than branched coverings of the sphere characterized by some finite data describing their singularities.
 
           In the comments to his paper in his \emph{Collected Works} Ahlfors writes: ``My aim was to interpret Nevanlinna's theorems as geometric properties of covering surfaces." He adds: ``My paper met with immediate recognition and earned me the Fields medal." In the presentation of Ahlfors' work that we already mentioned, Carathéodory declared that it is hard to say what is the most surprising thing among the following two: that Nevanlinna was able to develop his whole theory without any geometry, or that Ahlfors could concentrate the whole of Nevanlinna's theory in fourteen pages.

 We mention, as a later substantial step in the geometrization of value distribution\index{value distribution theory} theory, Chern's paper \cite{Chern}, in which value distribution is considered from a purely differential geometric point of view.  Chern studies  in this paper defect relations of complex analytic mappings between compact Riemann surfaces. He shows in particular that the number of points that are not covered by such a mapping is bounded above by the negative of the Euler characteristic of the range surface. This paved the way to a new series of results that relate value distribution to geometry and topology.

\section{Teichm\"uller's work}

The main lemma in Teichm\"uller's paper \cite{T13} gives conditions on a quasiconformal mapping $f(z)$ that insure that   $ \left| f(z)\right| \sim \textrm{const} \cdot  \left| z \right|$
as $z\to\infty$. Like the so-called Modulsatz, it gives a condition  on certain loci to be close to circles. In particular, these conditions allow the use of techniques of moduli of annuli in an efficient way in the inverse problem of value distribution theory. 

In fact, Teichm\"uller's contribution in \cite{T13}, as he put it himself, was motivated by Nevanlinna's inverse problem.\index{inverse problem (Nevanlinna)} The techniques he introduced in that paper were used by Goldberg  \cite[ch. 7]{Goldberg-Ostrovskii}, Le-Van Thiem \cite{Le, Le2}  and Drasin \cite{Drasin1977}   in their work on value distribution theory.

In his paper \cite{T33}, Teichm\"uller starts by recalling the uniformization map of a simply connected Riemann surfaces. He  then writes that the study of the value distribution of this unique function is one of the fundamental problems of modern function theory and that we are still quite far away from its solution.

Teichm\"uller restricts the problem to the case of Riemann surfaces defined as a ramified covering of the sphere. The value distribution of the uniformization map  of such surfaces has been investigated before, in a few cases, and he mentions the works on this problem by Nevanlinna \cite{N2},
Ahlfors \cite{Ahlfors-58},
Elfving \cite{Elfving-Uber},
Ullrich \cite{Ullrich-Jahr, Ullrich1939a},
and himself \cite{T9, T-Vermut}.

Surfaces with $p$ logarithmic ends were studied using line complexes by Nevanlinna, Ahlfors,   Elflving, Ullrich, Teichm\"uller and others. We survey this theory in another chapter of this volume, dedicated to the type problem \cite{ABP1}. Teichm\"uller considered the so-called \emph{surfaces with $p$ periodic ends}. These are ramified covers of the sphere with the property that outside a compact set (the ``kernel"), the ramification is periodic. These surfaces were first studied by Ullrich in \cite{Ullrich}. The line complex associated to such a surface has a certain periodic structure: there is kernel complex in which, at $p$ disjoint places, partial complexes are inserted, with this structure repeating periodically up to infinity. The surfaces with $p$ periodic ends that Teichm\"uller studies in the paper \cite{T33} are all of parabolic type.  He gives a method for constructing such surfaces and he studies their value distribution.\index{value distribution theory} Hans K\"unzi, in his two papers \cite{Kunzi-CR1}  and  \cite{Kunzi-CR2} (1952), gave another method for determining the characteristic, ramification indices and exceptional values for some surfaces studied by Teichm\"uller with a finite number of doubly periodic ends or of singly and doubly periodic ends.

The main result of Teichm\"uller's paper \cite{T33} is a set of  examples for which he gives explicit formulae for the counting functions $n(r,a)$, $N(r,a)$ and $m(r,a)$,  the characteristic function $T(r,f)$ and the deficiency function $\delta(a)$. The proofs use Nevanlinna's  two theorem. The theory of elliptic functions is used. Line complexes are used to construct the Riemann surface associated to the function $f$, and an explicit formula for the generating function associated to the Riemann surface obtained gives the information about the counting and the characteristic functions of $f$.

 Before the publication of his paper \cite{T33}, Teichm\"uller wrote the paper \cite{T8} (1937) in which he obtained a result on the generalized form of Nevanlinna's inverse problem.\index{inverse problem (Nevanlinna)!generalized}  Using the Koebe distortion theorem, he found conditions under which the proximity function $m(r,a)$ is bounded, and  the Valiron upper deficiency  $\Delta(a)$ is equal to $0$. His proof makes use of the Riemann surface associated to the function.  Teichm\"uller's theorem was generalized in 1948 by Collingwood in his paper \emph{Exceptional values of meromorphic functions}  \cite{Collingwood1949}, ``with only minor adaptations to the generalized conditions" (\cite{Collingwood1949} p. 311).    

 \section{Gauss--Lucas and Thurston}

           There are several natural questions concerning the location of the roots of a given polynomial. These questions are generally geometrical in nature, they are addressed by geometers, and they are not considered as belonging to the field of ``value distribution theory," because they concern the distribution of values of polynomials (and not of more general meromorphic functions). I would like to conclude this chapter by  mentioning one of these questions which is clearly related to that field.

A theorem attributed to Gauss\index{Gauss, Carl Friedrich (1777--1855 )} and Lucas\index{Lucas, Edouard (1842--1891)} says that the convex hull of the roots of a polynomial contains the roots of its derivative.\footnote{Gauss used implicitly this result in 1836 and Lucas published a proof in 1874. References are given in \cite{Marden}.} Obtaining more information on this convex hull and on related geometric loci was the subject of discussions on mathoverflow, in which Thurston\index{Thurston, William P. (1946--2012)} was involved. In a post written in 2011, participating to a discussion that concerns the intersection of the convex hull of level sets  $\{z\vert Q(z) = w\}$ of a polynomial $Q$, he writes: ``By chance, I've discussed this question a bit with Tan Lei; she made some nice movies of how the convex hulls of level sets vary with $w$. (Also, it's fun to look at their diagrams interactively manipulated in Mathematica). If I get my thoughts organized I'll post an answer." Thurston passed away the following year and he never had the chance to post the answer. In a tribute to Thurston, published three months after his death in the French electronic journal \emph{Images des math\'ematiques} dedicated to the popularization of  mathematics, Arnaud Ch\'eritat and Tan Lei published an article in which they present  Thurston's result  which gives a complete geometric picture of the situation. The result says the following: \emph{Let $P$ be a non-constant polynomial. Let $F$ be a half-plane bounded by a support line of the convex hull of the roots of the derivative $P'$ of $P$ and not containing this convex hull and let $c$ a root of $P'$ contained on this support line. Then, we can find a connected region contained in $F$ on which $P$ is bijective and whose interior is sent by $P$ on a plane with a slit along a ray starting at $P(c)$.} 
In 2015, Arnaud Ch\'eritat and Tan Lei published another article, in collaboration with Yan Gao, Yafei Ou,  titled ``A refinement of the Gauss--Lucas theorem (after W.P. Thurston)" \cite{Cheritat} in which they give more complete proofs of Thurston's theorem.\index{Gauss--Lucas theorem}\index{Theorem!Gauss--Lucas}

 \end{document}